# THE LONG-RUN BEHAVIOR OF THE STOCHASTIC REPLICATOR DYNAMICS

By Lorens A. Imhof

*Aachen University*

Fudenberg and Harris' stochastic version of the classical replicator dynamics is considered. The behavior of this diffusion process in the presence of an evolutionarily stable strategy is investigated. Moreover, extinction of dominated strategies and stochastic stability of strict Nash equilibria are studied. The general results are illustrated in connection with a discrete war of attrition. A persistence result for the maximum effort strategy is obtained and an explicit expression for the evolutionarily stable strategy is derived.

**1. Introduction.** The deterministic replicator dynamics is one of the most widely used dynamical models to describe the evolution of a population under selection. The evolution is governed by a symmetric two-player game with $n$ pure strategies, $1, \ldots, n$. Let $a_{jk}$ denote the pay-off to a player using strategy $j$ against an opponent playing strategy $k$. Let $A = (a_{jk})$. Suppose that every individual of the population is programmed to play one fixed pure strategy. For every point of time $t \geq 0$, let $\zeta_j(t)$ denote the size of the subpopulation whose individuals play strategy $j$, and let $\xi_j(t) = \zeta_j(t)/[\zeta_1(t) + \cdots + \zeta_n(t)]$ denote the proportion of $j$-players in the population. If the population state is $\xi(t) = (\xi_1(t), \ldots, \xi_n(t))^T$, then $\{A\xi(t)\}_j$ is the average pay-off to individuals playing $j$, when individuals are paired at random. Suppose that the pay-off represents the increase of fitness, measured as the number of offspring per unit of time. Then

$$\frac{d\zeta_j(t)}{dt} = \zeta_j(t)\{A\xi(t)\}_j, \qquad j = 1, \ldots, n, \tag{1.1}$$

and so

$$\frac{d\xi_j(t)}{dt} = \xi_j(t)[\{A\xi(t)\}_j - \xi(t)^T A\xi(t)], \qquad j = 1, \ldots, n. \tag{1.2}$$

This is the deterministic replicator dynamics of Taylor and Jonker (1978). See Hofbauer and Sigmund (1998) and Nowak and Sigmund (2004) for detailed discussions from a biological point of view and Weibull (1995) for a description in an economic context. See also Hofbauer and Sigmund (2003) for an extensive survey of deterministic evolutionary game dynamics.

Recently, models of evolutionary dynamics which incorporate stochastic effects have attracted substantial interest. The seminal paper of Foster and Young (1990) seems to be the first that presents a continuous-time replicator model based on a stochastic differential equation. Kandori, Mailath and Rob (1993) study a related discrete-time system. The present paper investigates the stochastic replicator dynamics introduced by Fudenberg and Harris (1992). This model is related to that of Foster and Young, but exhibits a boundary behavior that appears to be more realistic from a biological perspective. Following Fudenberg and Harris (1992), consider the stochastic variant of (1.1),

$$(1.3) \qquad dZ_j(t) = Z_j(t)[\{AX(t)\}_j \, dt + \sigma_j \, dW_j(t)], \qquad j = 1, \ldots, n,$$

where $(W_1(t), \ldots, W_n(t))^T = W(t)$ is an $n$-dimensional Brownian motion, $\sigma_1, \ldots, \sigma_n$ are positive coefficients and

$$X(t) = (X_1(t), \ldots, X_n(t))^T = \frac{1}{Z_1(t) + \cdots + Z_n(t)} (Z_1(t), \ldots, Z_n(t))^T.$$

The evolution of the population state $X(t)$ is then given by the stochastic replicator dynamics

$$(1.4) \qquad dX(t) = \mathbf{b}(X(t)) \, dt + C(X(t)) \, dW(t),$$

where

$$\mathbf{b}(\mathbf{x}) = [\mathrm{diag}(x_1, \ldots, x_n) - \mathbf{xx}^T][A - \mathrm{diag}(\sigma_1^2, \ldots, \sigma_n^2)]\mathbf{x}$$

and

$$C(\mathbf{x}) = [\mathrm{diag}(x_1, \ldots, x_n) - \mathbf{xx}^T] \mathrm{diag}(\sigma_1, \ldots, \sigma_n)$$

for $\mathbf{x} \in \Delta = \{\mathbf{y} \in (0,1)^n : y_1 + \cdots + y_n = 1\}$. In many interesting situations, the deterministic differential equation (1.2) has a stationary point in $\Delta$, which corresponds to a population state where every pure strategy is present. In fact every Nash equilibrium is stationary. On the other hand, the only stationary points for the stochastic differential equation (1.4) are the vertices of $\Delta$, corresponding to populations consisting of one common type of players.

A series of important results on the behavior of the stochastic replicator dynamics have been established for the case where the underlying game has two pure strategies. For example, Fudenberg and Harris (1992) and Saito (1997) examine properties of ergodic distributions, Amir and Berninghaus (1998) establish a result on equilibrium selection and Corradi and Sarin (2000) provide an asymptotic analysis. However, a large part of the arguments used there is tailored to the case $n = 2$ and cannot be extended to the general case $n > 2$. This is because when $n = 2$ one basically deals with one-dimensional diffusion processes, and many of the tools available for these processes are not applicable to higher-dimensional diffusions, which correspond to games with three or more pure strategies. In particular, in the general case, an approach via analyzing a closed form expression of the stationary distribution is not possible.

The present paper investigates (1.4) in the general case $n \geq 2$. Section 2 establishes a connection between stable behavior of the processes $X(t)$ and the static concept of an evolutionarily stable strategy (ESS), which has been introduced by Maynard Smith and Price (1973). Under suitable conditions, it is shown that if an ESS exists, then $X(t)$ is recurrent and the stationary distribution concentrates mass in a small neighborhood of the ESS. Explicit bounds for the expected time to reach that neighborhood are also given. Section 3 investigates dominated strategies. It is shown that the probability that the frequency of a dominated strategy is above a prescribed level decreases exponentially quickly to zero. Interestingly, it turns out that, depending on the sizes of the stochastic terms, weakly dominated strategies may become extinct in the stochastic model (1.4) even if they survive in the deterministic model (1.2). In Section 4 a sufficient condition is derived for a Nash equilibrium to be asymptotically stochastically stable. In this connection another example emerges which shows that the deterministic model and the stochastic model can lead to quite different predictions: In the Prisoner's Dilemma, the strategy "defect" is a strict Nash equilibrium and becomes predominant under (1.2), but may become extinct under (1.4).

By way of illustration, a discrete variant of the war of attrition is analyzed in some detail in the last section. This is a model which describes conflicts that are settled by display rather than violence; see Maynard Smith (1982). A rather general theorem on the persistence of the maximum effort strategy is obtained as a consequence of the results in Section 2. Furthermore, explicit expressions for ESSs are derived; the ESSs are given in terms of linear combinations of Chebyshev polynomials of the second kind evaluated along the imaginary axis. This yields a fairly accurate picture of the long-run behavior of the stochastic replicator dynamics when the conflicts are modeled by a war of attrition.

Hofbauer and Sigmund [(1998), Section 7.5] show that the deterministic replicator equation is, in a sense, equivalent to the deterministic Lotka–Volterra equation. The behavior of solutions to this equation under random perturbations has recently been investigated by Khasminskii and Klebaner (2001), Klebaner and Liptser (2001) and Skorokhod, Hoppensteadt and Salehi [(2002), Section 11]. There is almost no overlap with the results presented here.

**2. Stochastic replicator dynamics and evolutionarily stable strategies.**
The concept of a Nash equilibrium is too weak to yield reasonable convergence or stability results for (1.4). A somewhat stronger concept, which is of fundamental importance in evolutionary game theory, is that of an evolutionarily stable strategy (ESS), introduced by Maynard Smith and Price (1973). The closure $\overline{\Delta}$ of $\Delta$ is also referred to as the set of mixed strategies. A strategy $\mathbf{p} \in \overline{\Delta}$ is said to be an ESS if the following two conditions hold:

(i) $\mathbf{p}^T A \mathbf{p} \geq \mathbf{q}^T A \mathbf{p}$ for all $\mathbf{q} \in \overline{\Delta}$,

and

(ii) if $\mathbf{q} \neq \mathbf{p}$ and $\mathbf{p}^T A \mathbf{p} = \mathbf{q}^T A \mathbf{p}$, then $\mathbf{p}^T A \mathbf{q} > \mathbf{q}^T A \mathbf{q}$.

This static concept lies between that of a Nash and a strict Nash equilibrium, and turns out to be particularly relevant to the long-run analysis of (1.4).

For $\mathbf{x} \in \Delta$, let $P_{\mathbf{x}}$ denote the probability measure corresponding to the process $X(t)$ with initial condition $X(0) = \mathbf{x}$, and let $E_{\mathbf{x}}$ denote expectation with respect to $P_{\mathbf{x}}$. Note that $P_{\mathbf{x}}\{X(t) \in \Delta \text{ for all } t \geq 0\} = 1$ for all $\mathbf{x} \in \Delta$. Let $P(t, \mathbf{x}, G) = P_{\mathbf{x}}\{X(t) \in G\}$ for all Borel subsets $G \subset \Delta$. Let $\tau_G = \inf\{t > 0 : X(t) \in G\}$. For $\delta > 0$, let $U_\delta(\mathbf{x}) = \{\mathbf{y} \in \Delta : \|\mathbf{y} - \mathbf{x}\| < \delta\}$, where $\|\cdot\|$ denotes the Euclidean norm. Let $\mathbf{e}_j$ denote the $j$th unit vector in $\mathbb{R}^n$ and let $\mathbf{1} \in \mathbb{R}^n$ denote the vector all of whose entries are 1. The mixed strategy $\mathbf{e}_j$ is identified with the pure strategy $j$. The matrix $A$ is said to be conditionally negative definite if

$$\mathbf{y}^T A \mathbf{y} < 0 \qquad \text{for all } \mathbf{y} \in \mathbb{R}^n \text{ such that } \mathbf{1}^T \mathbf{y} = 0, \mathbf{y} \neq \mathbf{0}.$$

THEOREM 2.1. *Let $X(t)$ be given by the stochastic replicator dynamics* (1.4) *and let $\mathbf{p} \in \Delta$ be an ESS for the underlying pay-off matrix $A$. Set $\overline{A} = \frac{1}{2}(A + A^T)$, and let $\lambda_2$ be the second largest eigenvalue (counting multiplicity) of*

$$\overline{A} - \frac{1}{n}\overline{A}\mathbf{1}\mathbf{1}^T - \frac{1}{n}\mathbf{1}\mathbf{1}^T\overline{A} + \frac{\mathbf{1}^T A \mathbf{1}}{n^2}\mathbf{1}\mathbf{1}^T.$$

*Then*

(2.1) $$\lambda_2 < 0.$$

*Define $\kappa > 0$ by*

$$\kappa^2 = \frac{1}{2} \sum_{j=1}^{n} p_j \sigma_j^2 - \frac{1}{2 \sum_{j=1}^{n} \sigma_j^{-2}},$$

*and suppose that*

(2.2) $$\kappa < \frac{n}{n-1} \sqrt{|\lambda_2|} \min_{1 \leq j \leq n} p_j.$$

*Then $X(t)$ is recurrent, there exists a unique invariant probability measure $\pi$ on $\Delta$, and for every initial value $\mathbf{x} \in \Delta$, the transition probabilities $P(t, \mathbf{x}, \cdot)$ converge to $\pi$ in total variation. Moreover, for every $\delta > \kappa/\sqrt{|\lambda_2|}$,*

$$
\begin{aligned}
&\text{(a)} \quad \pi\{U_\delta(\mathbf{p})\} \geq 1 - \frac{\kappa^2}{|\lambda_2|\delta^2}, \\
&\text{(b)} \quad E_\mathbf{x} \tau_{\overline{U}_\delta(\mathbf{p})} \leq \frac{d(\mathbf{x}, \mathbf{p})}{|\lambda_2|\delta^2 - \kappa^2},
\end{aligned}
\tag{2.3}
$$

*and for every $t > 0$,*

$$
E_\mathbf{x} \frac{1}{t} \int_0^t \|X(s) - \mathbf{p}\|^2 \, ds \leq \frac{1}{|\lambda_2|} \left\{ \frac{d(\mathbf{x}, \mathbf{p})}{t} + \kappa^2 \right\},
\tag{2.4}
$$

*where $d(\mathbf{x}, \mathbf{p}) = \sum_{j\,:\,p_j > 0} p_j \log(p_j/x_j)$ is the Kullback–Leibler distance between $\mathbf{x}$ and $\mathbf{p}$.*

*Inequalities (2.1), (2.3)(b) and (2.4) also hold if the ESS $\mathbf{p} \in \overline{\Delta}$, provided that $A$ is conditionally negative definite.*

REMARK 2.1. The quantity $|\lambda_2|$ can be interpreted as a measure of how strongly the ESS $\mathbf{p}$ attracts $X(t)$ to a neighborhood of $\mathbf{p}$.

REMARK 2.2. Foster and Young (1990) point out that, in view of its local character, the ESS condition is not "quite the right concept of dynamical stability in a biological context." It is therefore not surprising that in the above theorem the ESS condition is augmented by some additional requirement: that $\kappa$ be not too large and that $A$ be conditionally negative definite if $\mathbf{p} \in \partial \Delta$. The second condition is easily seen to be satisfied in the examples in Section 5. Bapat and Raghavan [(1997), Section 4.1] provide some criteria to check whether a given matrix is conditionally negative definite.

The proof of Theorem 2.1 requires the following auxiliary result.

LEMMA 2.1. *Let $A \in \mathbb{R}^{n \times n}$, $n \geq 2$, be a conditionally negative definite matrix and let $\lambda_2$ be the second largest eigenvalue of*

$$
D := \overline{A} - \frac{1}{n}\overline{A}\mathbf{1}\mathbf{1}^T - \frac{1}{n}\mathbf{1}\mathbf{1}^T\overline{A} + \frac{\mathbf{1}^T A \mathbf{1}}{n^2}\mathbf{1}\mathbf{1}^T,
$$

*where $\overline{A} = \frac{1}{2}(A + A^T)$. Then*

$$
\max_{\substack{\mathbf{x}^T \mathbf{1} = 0 \\ \mathbf{x} \neq \mathbf{0}}} \frac{\mathbf{x}^T A \mathbf{x}}{\mathbf{x}^T \mathbf{x}} = \lambda_2 < 0.
$$

PROOF. Note first that

(2.5) $\quad \mathbf{x}^T D \mathbf{x} = \mathbf{x}^T A \mathbf{x} \quad$ for all $\mathbf{x} \in \mathbb{R}^n$ such that $\mathbf{1}^T \mathbf{x} = 0$.

The vector $\mathbf{1}$ is an eigenvector of $D$ corresponding to the eigenvalue $\lambda_1 = 0$. Thus if $\lambda$ is another eigenvalue of $D$ with corresponding eigenvector $\mathbf{y}$, then $\mathbf{1}^T \mathbf{y} = 0$. It then follows from (2.5) and the assumption that $A$ is conditionally negative definite that

$$\lambda = \frac{\mathbf{y}^T D \mathbf{y}}{\mathbf{y}^T \mathbf{y}} = \frac{\mathbf{y}^T A \mathbf{y}}{\mathbf{y}^T \mathbf{y}} \leq 0.$$

Thus $\lambda_1 = 0$ is the largest eigenvalue of $D$, and the variational description of $\lambda_2$, the second largest eigenvalue, and (2.5) yield

$$\lambda_2 = \max_{\substack{\mathbf{x}^T \mathbf{1} = 0 \\ \mathbf{x} \neq \mathbf{0}}} \frac{\mathbf{x}^T D \mathbf{x}}{\mathbf{x}^T \mathbf{x}} = \max_{\substack{\mathbf{x}^T \mathbf{1} = 0 \\ \mathbf{x} \neq \mathbf{0}}} \frac{\mathbf{x}^T A \mathbf{x}}{\mathbf{x}^T \mathbf{x}} < 0. \quad \square$$

PROOF OF THEOREM 2.1. Let $L$ denote the second-order differential operator associated with $X(t)$, that is,

(2.6) $\quad Lf(\mathbf{x}) = \sum_{j=1}^n b_j(\mathbf{x}) \frac{\partial f(\mathbf{x})}{\partial x_j} + \frac{1}{2} \sum_{j,k=1}^n \gamma_{jk}(\mathbf{x}) \frac{\partial^2 f(\mathbf{x})}{\partial x_j \, \partial x_k}, \quad f \in C^2(\Delta),$

where

$$b_j(\mathbf{x}) = x_j (\mathbf{e}_j - \mathbf{x})^T [A - \mathrm{diag}(\sigma_1^2, \ldots, \sigma_n^2)] \mathbf{x},$$

$$\gamma_{jk}(\mathbf{x}) = \sum_{\nu=1}^n c_{j\nu}(\mathbf{x}) c_{k\nu}(\mathbf{x}),$$

$$c_{jk}(\mathbf{x}) = \begin{cases} x_j(1 - x_j)\sigma_j, & j = k, \\ -x_j x_k \sigma_k, & j \neq k. \end{cases}$$

Suppose first that $\mathbf{p} \in \Delta$, and set $g(\mathbf{x}) = d(\mathbf{x}, \mathbf{p}) = \sum_j p_j \log(p_j/x_j)$ for all $\mathbf{x} \in \Delta$. Then, for all $\mathbf{x}$, $g(\mathbf{x}) \geq 0$ and

$$Lg(\mathbf{x}) = -\sum_{j=1}^n p_j (\mathbf{e}_j - \mathbf{x})^T [A - \mathrm{diag}(\sigma_1^2, \ldots, \sigma_n^2)] \mathbf{x}$$

$$+ \tfrac{1}{2} \sum_{j=1}^n p_j \left( \sigma_j^2 - 2 x_j \sigma_j^2 + \sum_{k=1}^n x_k^2 \sigma_k^2 \right)$$

$$= (\mathbf{x} - \mathbf{p})^T A \mathbf{x} - \tfrac{1}{2} \sum_{j=1}^n x_j^2 \sigma_j^2 + \tfrac{1}{2} \sum_{j=1}^n p_j \sigma_j^2.$$

As **p** is evolutionarily stable, $(\mathbf{x}-\mathbf{p})^T A\mathbf{p} \leq 0$. Since $\mathbf{p} \in \Delta$, $A$ is conditionally negative definite. This follows from the proof of Haigh's theorem (1975). Hence, in view of Lemma 2.1,

$$(\mathbf{x}-\mathbf{p})^T A\mathbf{x} \leq (\mathbf{x}-\mathbf{p})^T A(\mathbf{x}-\mathbf{p}) \leq \lambda_2 \|\mathbf{x}-\mathbf{p}\|^2$$

and $\lambda_2 < 0$. The Cauchy–Schwarz inequality gives $1 \leq (\sum_{j=1}^n x_j^2 \sigma_j^2) \sum_{j=1}^n \sigma_j^{-2}$, so that $-\sum_{j=1}^n x_j^2 \sigma_j^2 \leq -(\sum_{j=1}^n \sigma_j^{-2})^{-1}$. It now follows that

(2.7) $$Lg(\mathbf{x}) \leq \lambda_2 \|\mathbf{x}-\mathbf{p}\|^2 + \kappa^2, \qquad \mathbf{x} \in \Delta.$$

Suppose that $\delta^2 > \kappa^2/|\lambda_2|$. For every $\mathbf{x} \in \Delta \setminus U_\delta(\mathbf{p})$, $Lg(\mathbf{x}) \leq \lambda_2 \delta^2 + \kappa^2$, and it follows by Itô's formula that $g(X(t)) - (\lambda_2 \delta^2 + \kappa^2)t$ is a local supermartingale on $[0, \tau_{\overline{U}_\delta(\mathbf{p})})$. Hence [cf. proof of Theorem 5.3 in Durrett (1996), page 268] $g(\mathbf{x}) \geq (|\lambda_2|\delta^2 - \kappa^2)E_\mathbf{x} \tau_{\overline{U}_\delta(\mathbf{p})}$, proving (2.3)(b).

To prove recurrence, consider the transformed process $Y(t) = \Psi(X(t))$, where $\Psi: \Delta \to \mathbb{R}^{n-1}$ is defined by $\Psi(\mathbf{x}) = (\log(x_1/x_n), \ldots, \log(x_{n-1}/x_n))^T$. One has

$$dY_j(t) = \{(\mathbf{e}_j - \mathbf{e}_n)^T A \Psi^{-1}(Y(t)) - \tfrac{1}{2}(\sigma_j^2 - \sigma_n^2)\} dt + \sigma_j dW_j(t) - \sigma_n dW_n(t),$$
$$j = 1, \ldots, n-1,$$

where $\Psi^{-1}(\mathbf{y}) = (1 + e^{y_1} + \cdots + e^{y_{n-1}})^{-1}(e^{y_1}, \ldots, e^{y_{n-1}}, 1)^T$. Note that the second-order differential operator associated with $Y(t)$ is uniformly elliptic. It will next be shown that $\partial U_{\delta_0}$ has positive distance from $\partial \Delta$ for some $\delta_0 > \kappa/\sqrt{|\lambda_2|}$. This implies that $\Psi(U_{\delta_0})$ is a compact set. In view of (2.3)(b) it will then follow that $Y(t)$ is recurrent [Bhattacharya (1978) and Khas'minskii (1960)] and the transition probabilities converge in total variation to the unique stationary probability measure [Durrett (1996), Chapter 7]. The same applies then to $X(t)$. By (2.2), one may choose $\delta_0$ such that $\kappa/\sqrt{|\lambda_2|} < \delta_0 < (n/(n-1))\min_j p_j$. Suppose $\mathbf{y} \in \mathbb{R}^n$, $\sum_{j=1}^n y_j = 1$ and $\|\mathbf{y}-\mathbf{p}\| = \delta_0$. Let $j_0$ be such that $|y_{j_0} - p_{j_0}| = \max_{1 \leq j \leq n} |y_j - p_j|$, and set

$$\mathbf{z} = \left(y_{j_0} - p_{j_0}, \frac{p_{j_0} - y_{j_0}}{n-1}, \ldots, \frac{p_{j_0} - y_{j_0}}{n-1}\right)^T \in \mathbb{R}^n.$$

One may verify that $\mathbf{z}$ is majorized by $\mathbf{y} - \mathbf{p}$ in the sense of Definition A.1 in Marshall and Olkin [(1979), page 7], and it follows from Proposition C.1 in Marshall and Olkin [(1979), page 64] that

$$\delta_0^2 = \sum_{j=1}^n (y_j - p_j)^2 \geq \sum_{j=1}^n z_j^2 = \left(1 + \frac{1}{n-1}\right)(y_{j_0} - p_{j_0})^2$$
$$= \frac{n}{n-1} \max_{1 \leq j \leq n} (y_j - p_j)^2.$$

Thus, for every $j$, $y_j \geq \min_{1 \leq k \leq n} p_k - (n-1)\delta_0/n > 0$, showing that the distance between $\partial U_{\delta_0}$ and $\partial \Delta$ is positive.

For $K > g(\mathbf{x})$ let $\tilde{\tau}_K = \inf\{t > 0 : g(X(t)) = K\}$. Then, by Dynkin's formula and (2.7),

$$0 \leq E_{\mathbf{x}} g(X(t \wedge \tilde{\tau}_K)) = g(\mathbf{x}) + E_{\mathbf{x}} \int_0^{t \wedge \tilde{\tau}_K} Lg(X(s)) \, ds$$

$$\leq g(\mathbf{x}) + \lambda_2 E_{\mathbf{x}} \int_0^{t \wedge \tilde{\tau}_K} \|X(s) - \mathbf{p}\|^2 \, ds + \kappa^2 E_{\mathbf{x}}(t \wedge \tilde{\tau}_K).$$

If $K \to \infty$, then $t \wedge \tilde{\tau}_k \to t$, and (2.4) follows by the bounded convergence theorem. To prove (2.3)(a), let $\chi_{\overline{U}_\delta^C(\mathbf{p})}$ denote the indicator function of $\overline{U}_\delta^C(\mathbf{p}) = \Delta \setminus \overline{U}_\delta(\mathbf{p})$. Then, by (2.4),

$$\pi(\overline{U}_\delta^C(\mathbf{p})) = \lim_{t \to \infty} E_{\mathbf{x}} \frac{1}{t} \int_0^t \chi_{\overline{U}_\delta^C(\mathbf{p})}(X(s)) \, ds$$

$$\leq \lim_{t \to \infty} E_{\mathbf{x}} \frac{1}{t} \int_0^t \frac{\|X(s) - \mathbf{p}\|^2}{\delta^2} \, ds \leq \frac{\kappa^2}{|\lambda_2|\delta^2}.$$

An inspection of the above arguments shows that if $A$ is conditionally negative definite, then (2.1), (2.3)(b) and (2.4) hold if the ESS $\mathbf{p} \in \partial \Delta$. In this case, however, $X(t)$ need not be recurrent. $\square$

If $\mathbf{p}$ is an ESS for $A$, then there exists a constant $c \in \mathbb{R}$ such that $\{A\mathbf{p}\}_j = c$ for all $j \in \{1, \ldots, n\}$ with $p_j > 0$; see Hofbauer and Sigmund [(1998), page 63]. Thus

$$[\text{diag}(p_1, \ldots, p_n) - \mathbf{p}\mathbf{p}^T]A\mathbf{p} = \mathbf{0},$$

so that the drift vector $\mathbf{b}(\mathbf{x})$ of the stochastic differential equation (1.4) will in general not be zero at $\mathbf{x} = \mathbf{p}$. If, however, $\mathbf{p}$ is an ESS for the modified pay-off matrix $B = A - \text{diag}(\sigma_1^2, \ldots, \sigma_n^2)$, then $\mathbf{b}(\mathbf{p}) = \mathbf{0}$. From this point of view it is more natural to investigate the distance between $X(t)$ and an ESS for $B$. A simple modification of the proof of Theorem 2.1 yields the following result. The analogous results on recurrence and the stationary distribution are omitted for brevity.

THEOREM 2.2. *Let $X(t)$ be given by the stochastic replicator dynamics (1.4) with underlying pay-off matrix $A$. Let $\mathbf{p} \in \overline{\Delta}$ be an ESS for the modified pay-off matrix $A - \text{diag}(\sigma_1^2, \ldots, \sigma_n^2)$. Suppose also that $A - \text{diag}(\frac{1}{2}\sigma_1^2, \ldots, \frac{1}{2}\sigma_n^2)$ is conditionally negative definite. Then for every initial state $\mathbf{x} \in \Delta$ and every $t > 0$,*

$$(2.8) \quad E_{\mathbf{x}} \frac{1}{t} \int_0^t \|X(s) - \mathbf{p}\|^2 \, ds \leq \frac{1}{|\lambda_2'|} \left\{ \frac{d(\mathbf{x}, \mathbf{p})}{t} + \frac{1}{2} \sum_{j=1}^n p_j(1 - p_j)\sigma_j^2 \right\},$$

where $d(\mathbf{x}, \mathbf{p}) = \sum_{j : p_j > 0} p_j \log(p_j/x_j)$ and $\lambda_2'$ is the second largest eigenvalue of

$$\overline{A} - \frac{1}{n}\overline{A}\mathbf{1}\mathbf{1}^T - \frac{1}{n}\mathbf{1}\mathbf{1}^T\overline{A} + \frac{\mathbf{1}^T\overline{A}\mathbf{1}}{n^2}\mathbf{1}\mathbf{1}^T,$$

where $\overline{A} = \frac{1}{2}[A + A^T - \mathrm{diag}(\sigma_1^2, \ldots, \sigma_n^2)]$.

REMARK 2.3. If the ESS $\mathbf{p} \in \Delta$, it follows that $A - \mathrm{diag}(\sigma_1^2, \ldots, \sigma_n^2)$ is conditionally negative definite, so that in this case, the assumption that $A - \mathrm{diag}(\frac{1}{2}\sigma_1^2, \ldots, \frac{1}{2}\sigma_n^2)$ should be conditionally negative definite is not very restrictive. To compare (2.4) and (2.8) note that $|\lambda_2'| > |\lambda_2|$ and

$$\frac{1}{2}\sum_{j=1}^n p_j(1-p_j)\sigma_j^2 \leq -\frac{1}{2\sum_{j=1}^n \sigma_j^{-2}} + \frac{1}{2}\sum_{j=1}^n p_j\sigma_j^2.$$

**3. Extinction of dominated strategies.** This section is concerned with the evolution of strategies that are inferior to other strategies in the sense of domination.

A strategy $\mathbf{p} \in \overline{\Delta}$ is said to be weakly dominated by strategy $\mathbf{q} \in \overline{\Delta}$ if

$$\mathbf{p}^T A \mathbf{r} \leq \mathbf{q}^T A \mathbf{r} \qquad \text{for all } \mathbf{r} \in \overline{\Delta}$$

with strict inequality for some $\mathbf{r}$. If the inequality is strict for all $\mathbf{r}$, then $\mathbf{p}$ is said to be strictly dominated by $\mathbf{q}$.

For the deterministic replicator dynamics (1.2), Akin (1980) has shown that strictly dominated pure strategies become extinct; more precisely, their frequencies in the population converge to zero. Theorem 3.1 establishes that under the stochastic replicator dynamics (1.4) even pure strategies that are only weakly dominated become extinct under a suitable condition on the diffusion coefficients $\sigma_1, \ldots, \sigma_n$. Theorem 3.1 also gives an upper bound for the probability that at a given point of time $t$ the frequency of a dominated strategy is above a prescribed value $\varepsilon > 0$. The bound converges exponentially quickly to zero as $t \to \infty$.

THEOREM 3.1. *Let $X(t)$ be given by* (1.4). *Suppose that the pure strategy $k$ is weakly dominated by some mixed strategy* $\mathbf{p} \in \overline{\Delta}$. *Set* $c_1 = \min_{\mathbf{q} \in \overline{\Delta}} \mathbf{p}^T A \mathbf{q} - \mathbf{e}_k^T A \mathbf{q}$ *and suppose that $\sigma_1, \ldots, \sigma_n$ are such that*

(3.1) $$c_2 = -\frac{\sigma_k^2}{2} + \frac{1}{2}\sum_{j=1}^n p_j \sigma_j^2 < c_1.$$

*Then for every initial state* $\mathbf{x} \in \Delta$,

$$P_{\mathbf{x}}\{X_k(t) = o(\exp[-(c_1 - c_2)t + 3\sigma_{\max}\sqrt{t \log \log t}])\} = 1,$$

*and for $0 < \varepsilon < 1$ and $t > 0$,*

$$P_{\mathbf{x}}\{X_k(t) > \varepsilon\} < 1 - \Phi\left\{\frac{c_3(\mathbf{x}) + \log \varepsilon + (c_1 - c_2)t}{\sigma_{\max}\sqrt{2t}}\right\},$$

*where $\Phi(v)$ is the normal distribution function, $\sigma_{\max} = \max\{\sigma_1, \ldots, \sigma_n\}$ and $c_3(\mathbf{x}) = \sum_{j=1}^{n} p_j \log(x_j/x_k)$.*

PROOF. Let $H(t) = \log X_k(t) - \sum_{j=1}^{n} p_j \log X_j(t)$ for $t \geq 0$. Then, by Itô's formula,

$$H(t) = H(0) + \int_0^t \mathbf{e}_k^T A X(s) - \mathbf{p}^T A X(s) - \frac{\sigma_k^2}{2} + \frac{1}{2}\sum_{j=1}^{n} p_j \sigma_j^2 \, ds$$

$$+ \sigma_k W_k(t) - \sum_{j=1}^{n} p_j \sigma_j W_j(t)$$

$$\leq H(0) + (c_2 - c_1)t + \tilde{\sigma}\widetilde{W}(t),$$

where $\tilde{\sigma} = [(1-p_k)^2 \sigma_k^2 + \sum_{j \neq k} p_j^2 \sigma_j^2]^{1/2}$ and $\widetilde{W}(t) = [\sigma_k W_k(t) - \sum_{j=1}^{n} p_j \sigma_j \times W_j(t)]/\tilde{\sigma}$ is a standard Brownian motion. Clearly, $\tilde{\sigma} \leq \sqrt{2}\sigma_{\max}$. It follows that $P_{\mathbf{x}}$-almost surely,

$$\limsup_{t \to \infty} X_k(t) \exp[(c_1 - c_2)t - 3\sigma_{\max}\sqrt{t \log \log t}]$$

$$\leq \limsup_{t \to \infty} \exp[H(t) + (c_1 - c_2)t - 3\sigma_{\max}\sqrt{t \log \log t}]$$

$$\leq \limsup_{t \to \infty} \exp[-c_3(\mathbf{x}) + \tilde{\sigma}\widetilde{W}(t) - 3\sigma_{\max}\sqrt{t \log \log t}] = 0$$

by the law of the iterated logarithm. Moreover,

$$P_{\mathbf{x}}\{X_k(t) > \varepsilon\} \leq P_{\mathbf{x}}\{H(t) > \log \varepsilon\}$$

$$\leq P_{\mathbf{x}}\{-c_3(\mathbf{x}) + (c_2 - c_1)t + \tilde{\sigma}W(t) > \log \varepsilon\}$$

$$< 1 - \Phi\left\{\frac{c_3(\mathbf{x}) + \log \varepsilon + (c_1 - c_2)t}{\sigma_{\max}\sqrt{2t}}\right\}. \qquad \square$$

REMARK 3.1. Condition (3.1) is always satisfied if $k$ is strictly dominated by $\mathbf{p}$ and $\sigma_1 = \cdots = \sigma_n$. The condition is also satisfied when $k$ is merely weakly dominated and $\sigma_k > \sigma_j$ for every $j \neq k$. Thus if the diffusion coefficient corresponding to a weakly dominated strategy is large enough, its frequency converges to zero. This behavior is different from the behavior of weakly dominated strategies under the deterministic replicator dynamics where weakly dominated strategies may well persist with any prescribed positive population share; see Weibull [(1995), Example 3.4, page 84]. This

STOCHASTIC REPLICATOR DYNAMICS 11difference between the deterministic and the stochastic population dynamics agrees with the findings of Alvarez (2000) which show that, under mild conditions, "increased stochastic fluctuations decrease the expected population density." Cabrales (2000) proves that in a similar stochastic model iteratively strictly dominated strategies become rare, provided stochastic effects are sufficiently small.

REMARK 3.2. In the situation of Theorem 3.1,
$$P_{\mathbf{x}}\{X_k(t) > \varepsilon\} = o(e^{-\gamma t}), \qquad t \to \infty,$$
for any $0 < \gamma < (c_1 - c_2)^2/(4\sigma_{\max}^2)$. This is easily verified using the bound $1 - \Phi(v) \leq \exp(-v^2/2)$, $v > 0$.

REMARK 3.3. The assumptions of Theorem 3.1 may be satisfied even if the pure strategy $k$ is not dominated by any other pure strategy. For example, let

$$A = \begin{pmatrix} 2 & 2 & 2 \\ 4 & 1 & 1 \\ 1 & 4 & 4 \end{pmatrix} \quad \text{and} \quad \frac{\sigma_2^2 + \sigma_3^2}{2} < 1 + \sigma_1^2.$$

Then the pure strategy 1 is strictly dominated by $\mathbf{p} = (0, \frac{1}{2}, \frac{1}{2})^T$ and, according to Theorem 3.1, $X_1(t) \to 0$ almost surely, even though neither strategy 2 nor strategy 3 dominates strategy 1.

**4. Stochastic stability of Nash equilibria.** The last section dealt with the extinction of pure strategies that were inferior to at least one strategy. The present section investigates strategies that can be regarded as being locally superior to all other strategies. The relevant concept is that of a strict Nash equilibrium.

A strategy $\mathbf{p} \in \overline{\Delta}$ is called a Nash equilibrium if
$$\mathbf{p}^T A \mathbf{p} \geq \mathbf{q}^T A \mathbf{p} \qquad \text{for all } \mathbf{q} \in \overline{\Delta}.$$
If the inequality is strict for all $\mathbf{q} \neq \mathbf{p}$, then $\mathbf{p}$ is a strict Nash equilibrium.

In other words, a Nash equilibrium is a best reply to itself, and a strict Nash equilibrium is the unique best reply to itself. Only pure strategies can be strict Nash equilibria; see Weibull [(1995), page 15]. Thus if nearly the whole population plays a strict Nash equilibrium, then the highest pay-off is obtained by exactly that strategy so that natural selection would not favor any other strategy. This suggests that a strict Nash equilibrium should be an asymptotically stable state, which is indeed the case under the deterministic replicator dynamics (1.2). This need not be the case for the stochastic replicator dynamics (1.4), as is illustrated by the following example.

EXAMPLE. Consider the Prisoner's Dilemma game [Hofbauer and Sigmund (1998), page 101] with two pure strategies: $1 =$ co-operate, $2 =$ defect, and pay-off matrix

$$A = \begin{pmatrix} a_{11} & a_{12} \\ a_{21} & a_{22} \end{pmatrix}, \qquad a_{21} > a_{11} > a_{22} > a_{12}.$$

Here strategy 2 is a strict Nash equilibrium and under (1.2), $\lim_{t\to\infty} \xi_2(t) = 1$ for all initial states $\xi(0) \in \Delta$. On the other hand, if

$$\frac{\sigma_2^2}{2} > \frac{\sigma_1^2}{2} + \max\{a_{21} - a_{11}, a_{22} - a_{12}\},$$

then condition (3.1) of Theorem 3.1 is satisfied with $k = 2$, $\mathbf{p} = (1,0)^T$, so that $P_\mathbf{x}\{\lim_{t\to\infty} X_1(t) = 1\}$ for all $\mathbf{x} \in \Delta$. This is a reasonable behavior, because if all players co-operate, then the received pay-off, $a_{11}$, is larger than $a_{22}$, the pay-off they receive when all players defect. Thus the stochastic model may explain the spread of co-operative behavior, which could not be observed in the deterministic model. This fact agrees with results of Nowak, Sasaki, Taylor and Fudenberg (2004) who study a discrete-time Markov chain to explain the emergence of co-operation in a finite population that plays the Prisoner's Dilemma game.

The next theorem gives a sufficient condition for a strict Nash equilibrium to be asymptotically stochastically stable. Notice that a pure strategy $k$ is a strict Nash equilibrium if and only if $a_{kk} > a_{jk}$ for all $j \neq k$.

THEOREM 4.1. *Let $X(t)$ be given by (1.4). Let $k$ be a strict Nash equilibrium. Suppose that $\sigma_k$ is so small that*

(4.1) $$a_{kk} > a_{jk} + \sigma_k^2 \qquad \text{for all } j \neq k.$$

*Then $\mathbf{e_k}$ is asymptotically stochastically stable. That is, for any neighborhood $U$ of $\mathbf{e}_k$ and for any $\varepsilon > 0$ there is a neighborhood $V$ of $\mathbf{e_k}$ such that*

$$P_\mathbf{x}\left\{X(t) \in U \text{ for all } t \geq 0, \lim_{t\to\infty} X(t) = \mathbf{e}_k\right\} \geq 1 - \varepsilon$$

*for every initial state $\mathbf{x} \in V \cap \Delta$.*

PROOF. The proof is an application of the stochastic Lyapunov method. Consider the Lyapunov function $\phi(\mathbf{y}) = 1 - y_k$. Evidently, $\phi(\mathbf{y}) \geq 0$ for all $\mathbf{y} \in \overline{\Delta}$ with equality if and only if $\mathbf{y} = \mathbf{e}_k$. It will be shown that there is a constant $c > 0$ and a neighborhood $V_0$ of $\mathbf{e}_k$ such that

(4.2) $$L\phi(\mathbf{y}) \leq -c\phi(\mathbf{y}) \qquad \text{for all } \mathbf{y} \in V_0 \cap \Delta,$$

where $L$ is the differential operator given by (2.6). The assertion then follows from Theorem 4 and Remark 2 in Gichman and Skorochod [(1971), pages 314 and 315]. Write $B = A - \text{diag}(\sigma_1^2, \ldots, \sigma_n^2)$. For all $\mathbf{y} \in \Delta$,

$$\begin{aligned} L\phi(\mathbf{y}) &= -y_k(\mathbf{e}_k - \mathbf{y})^T B \mathbf{y} \\ &= y_k \sum_{\substack{\mu \neq k \\ \nu \neq k}} y_\mu b_{\mu\nu} y_\nu - y_k(1 - y_k) \sum_{\nu \neq k} b_{k\nu} y_\nu \\ &\quad + y_k^2 \bigg\{ -(1 - y_k) b_{kk} + \sum_{\mu \neq k} y_\mu b_{\mu k} \bigg\}. \end{aligned}$$

Let $\beta = \max\{|b_{\mu\nu}| : \mu, \nu = 1, \ldots, n\}$. Then

$$\sum_{\substack{\mu \neq k \\ \nu \neq k}} y_\mu b_{\mu\nu} y_\nu \leq \beta \sum_{\mu \neq k} y_\mu \sum_{\nu \neq k} y_\nu = \beta(1 - y_k)^2, \qquad -\sum_{\nu \neq k} b_{k\nu} y_\nu \leq \beta(1 - y_k).$$

Moreover, condition (4.1) ensures that for some $\alpha > 0$, $b_{\mu k} \leq b_{kk} - \alpha$ for all $\mu \neq k$, so that

$$-(1 - y_k) b_{kk} + \sum_{\mu \neq k} y_\mu b_{\mu k} \leq -(1 - y_k) b_{kk} + (b_{kk} - \alpha) \sum_{\mu \neq k} y_\mu = -\alpha(1 - y_k).$$

Hence

$$L\phi(\mathbf{y}) \leq 2\beta y_k(1 - y_k)^2 - \alpha y_k^2(1 - y_k) = -y_k\{(\alpha + 2\beta)y_k - 2\beta\}\phi(\mathbf{y}),$$

which proves (4.2) with $V_0 = \{\mathbf{y} \in \Delta : y_k > \frac{1}{2} \frac{\alpha + 4\beta}{\alpha + 2\beta}\}$ and $c = \frac{\alpha}{4}$. $\square$

REMARK 4.1. Condition (4.1) means that $k$ is a strict Nash equilibrium with respect to the modified pay-off matrix $B$. If $k$ is only a neutrally stable strategy with respect to $B$ [see Weibull (1995), Definition 2.4, page 46], then, by Proposition 2.7 of Weibull [(1995), page 48], $L\phi(\mathbf{y}) = -y_k(\mathbf{e}_k - \mathbf{y})^T B \mathbf{y} \leq 0$ for all $\mathbf{y}$ in a certain neighborhood of $\mathbf{e}_k$. Hence in this case Theorem 4 of Gichman and Skorochod [(1971), page 314] yields that $\mathbf{e}_k$ is still stochastically stable.

Theorem 4.1 says that if the population is in a state sufficiently near to a strict Nash equilibrium $\mathbf{e}_k$, then, with probability close to 1, that equilibrium will actually be selected by the stochastic replicator dynamics in the sense that $\lim_{t \to \infty} X(t) = \mathbf{e}_k$. If there are several strict Nash equilibria and the initial state is not close to any of them, it is neither clear which one will be selected nor in fact if any will be selected at all. The next theorem establishes that when the underlying game is a coordination game, that is, every pure stategy is a strict Nash equilibrium, then it is almost certain that one of the equilibria will be selected.

THEOREM 4.2. *Let $A$ be the pay-off matrix of a co-ordination game and let $X(t)$ be given by (1.4). Suppose that, for every $k$, $\sigma_k$ is so small that $a_{kk} > a_{jk} + \sigma_k^2$ for all $j \neq k$. Then, for every initial state $\mathbf{x} \in \Delta$,*

$$P_{\mathbf{x}}\left\{\lim_{t \to \infty} X(t) = \mathbf{e}_k \text{ for some } k\right\} = 1.$$

The proof hinges on the following theorem, which is of interest in its own right. It states that for any underlying game, $X(t)$ will come arbitrarily close to one of the points $\mathbf{e}_1, \ldots, \mathbf{e}_n$ in finite time.

THEOREM 4.3. *Let $A$ be an arbitrary pay-off matrix, let $X(t)$ be given by (1.4) and let $\mathbf{x} \in \Delta$. Let $\varepsilon > 0$. Consider the hitting time*

$$\tau_\varepsilon = \inf\{t > 0 : X_k(t) \geq 1 - \varepsilon \text{ for some } k \in \{1, \ldots, n\}\}.$$

*Then $E_{\mathbf{x}}\tau_\varepsilon < \infty$. Moreover,*

$$P_{\mathbf{x}}\left\{\sup_{t > 0} \max\{X_1(t), \ldots, X_n(t)\} = 1\right\} = 1.$$

PROOF. For $\alpha > 0$ and $\mathbf{y} \in \overline{\Delta}$ define

$$\phi_\alpha(\mathbf{y}) = \phi(\mathbf{y}) = ne^\alpha - \sum_{k=1}^n e^{\alpha y_k}.$$

Let $B = A - \text{diag}(\sigma_1^2, \ldots, \sigma_n^2)$ and let $L$ be given by (2.6). Then

$$L\phi(\mathbf{y}) = -\alpha \sum_{k=1}^n y_k (\mathbf{e}_k - \mathbf{y})^T B \mathbf{y} e^{\alpha y_k} - \frac{\alpha^2}{2} \sum_{k=1}^n y_k^2 \left\{\sigma_k^2 (1 - y_k)^2 + \sum_{j \neq k} \sigma_j^2 y_j^2\right\} e^{\alpha y_k}.$$

Let $\beta > 0$ be such that

$$|(\mathbf{e}_k - \mathbf{y})^T B \mathbf{y}| \leq \beta \qquad \text{for all } \mathbf{y} \in \Delta \text{ and } k \in \{1, \ldots, n\}.$$

Let $\sigma_{\min} = \min\{\sigma_1, \ldots, \sigma_n\}$. Then

$$L\phi(\mathbf{y}) \leq \alpha \sum_{k=1}^n y_k e^{\alpha y_k} \left\{\beta - \frac{\alpha \sigma_{\min}^2}{2} y_k (1 - y_k)^2\right\}.$$

Let $\alpha > 0$ be so large that

$$\alpha \frac{\sigma_{\min}^2}{2} y(1-y)^2 \geq n\beta + 1 \qquad \text{for all } y \in \left[\frac{1}{n}, 1 - \varepsilon\right].$$

Suppose that $\mathbf{y} \in \Delta$ is such that $y_k \leq 1 - \varepsilon$ for all $k = 1, \ldots, n$. Then there is at least one $y_k$ in $[\frac{1}{n}, 1 - \varepsilon]$, and so

$$L\phi(\mathbf{y}) \leq \alpha\beta \sum_{k: y_k < 1/n} y_k e^{\alpha y_k} + \alpha \sum_{k: y_k \in [1/n, 1-\varepsilon]} y_k e^{\alpha y_k}\{-(n-1)\beta - 1\}$$

$$\leq \alpha\beta(n-1)\frac{e^{\alpha/n}}{n} + \alpha\frac{e^{\alpha/n}}{n}\{-(n-1)\beta - 1\}$$

$$\leq -\alpha\frac{e^{\alpha/n}}{n}.$$

Thus, by Dynkin's formula, for every $T < \infty$,

$$0 \leq E_\mathbf{x}\phi\{X(\tau_\varepsilon \wedge T)\} = \phi(\mathbf{x}) + E_\mathbf{x}\int_0^{\tau_\varepsilon \wedge T} L\phi(X(s))\,ds$$

$$\leq ne^\alpha - \alpha\frac{e^{\alpha/n}}{n}E_\mathbf{x}(\tau_\varepsilon \wedge T).$$

Letting $T \to \infty$, one obtains by monotone convergence that $E_\mathbf{x}\tau_\varepsilon < n^2 e^\alpha/\alpha < \infty$.

Choosing $\varepsilon = 1/m$, one obtains in particular that

$$P_\mathbf{x}\left\{\sup_{t>0}\max\{X_1(t), \ldots, X_n(t)\} \geq 1 - \frac{1}{m}\right\} = 1$$

for every $m \in \mathbb{N}$. Hence

$$P_\mathbf{x}\left\{\sup_{t>0}\max\{X_1(t), \ldots, X_n(t)\} = 1\right\}$$

$$= P_\mathbf{x}\left\{\bigcap_{m=1}^\infty \left\{\sup_{t>0}\max\{X_1(t), \ldots, X_n(t)\} \geq 1 - \frac{1}{m}\right\}\right\} = 1. \quad \square$$

PROOF OF THEOREM 4.2. Let $\varepsilon > 0$ and suppose that $a_{kk} > a_{jk} + \sigma_k^2$ for all $j, k = 1, \ldots, n$ with $j \neq k$. Then, for every $k$, there exists, by Theorem 4.1, some $\delta_k > 0$ such that

$$P_\mathbf{x}\left\{\lim_{t\to\infty} X(t) = \mathbf{e}_k\right\} > 1 - \varepsilon \quad \text{if } x_k \geq 1 - \delta_k.$$

Set $\tau = \inf\{t \geq 0 : X_k(t) \geq 1 - \delta_k \text{ for some } k\}$ and $F = \{\lim_{t\to\infty} X(t) = \mathbf{e}_k \text{ for some } k\}$. Let $\chi_F$ denote the indicator function of $F$. According to Theorem 4.3, $\tau$ is $P_\mathbf{x}$-almost surely finite, and so, by the strong Markov property of Itô diffusions,

$$P_\mathbf{x}(F) = E_\mathbf{x}E_{X(\tau)}\chi_F \geq 1 - \varepsilon.$$

As $\varepsilon > 0$ was arbitrary, the assertion follows. $\square$

**5. A discrete war of attrition.** Theorems 2.1 and 2.2 show that an ESS, if it exists, gives precise information about the long-run behavior of the stochastic replicator dynamics. In this section explicit expressions are derived for ESSs for a discrete variant of the war of attrition, introduced by Maynard Smith and Price (1973). See Maynard Smith (1982) and Bishop and Cannings (1978) for a detailed discussion and extensions.

In the discrete symmetric war of attrition each player selects a pure strategy $j \in \{0, 1, \ldots, n\}$, which determines the maximum length of time, $c_j$, the player is willing to display for. The contest progresses until one of the players has reached his chosen limit; this player leaves and the other player obtains a reward. The value of the reward is constant or a decreasing function of the length of the contest. Both players incur a cost given by the length of the contest. If both players have chosen the same length of time, the reward is shared.

Specifically, the pay-off matrix $A = (a_{jk})$, $j, k = 0, \ldots, n$, is given by

$$(5.1) \quad a_{jk} = \begin{cases} v_k - c_k, & j > k, \\ \dfrac{v_k}{2} - c_k, & j = k, \\ -c_j, & j < k, \end{cases}$$

where

$$0 \leq c_0 < c_1 < \cdots < c_n \quad \text{and} \quad v_0 \geq v_1 \geq \cdots \geq v_n > 0.$$

The corresponding stochastic replicator dynamics is

$$(5.2) \quad dX(t) = \mathbf{b}(X(t)) \, dt + C(X(t)) \, dW(t),$$

where $W(t)$ denotes an $(n+1)$-dimensional Brownian motion and, for $\mathbf{x} = (x_0, \ldots, x_n)^T$,

$$\mathbf{b}(\mathbf{x}) = [\operatorname{diag}(x_0, \ldots, x_n) - \mathbf{x}\mathbf{x}^T][A - \operatorname{diag}(\sigma_0^2, \ldots, \sigma_n^2)]\mathbf{x}$$

and

$$C(\mathbf{x}) = [\operatorname{diag}(x_0, \ldots, x_n) - \mathbf{x}\mathbf{x}^T]\operatorname{diag}(\sigma_0, \ldots, \sigma_n).$$

Theorem 2.2 suggests to consider ESSs not only of $A$ but also of the modified pay-off matrix $B = (b_{jk})$, $j, k = 0, \ldots, n$, with

$$(5.3) \quad b_{jk} = \begin{cases} v_k - c_k, & j > k, \\ \dfrac{v_k}{2} - c_k - \rho_k, & j = k, \\ -c_j, & j < k, \end{cases}$$

and

$$0 \leq c_0 < c_1 < \cdots < c_n, \qquad v_0 \geq v_1 \geq \cdots \geq v_n > 0, \qquad 0 \leq \rho_k < \dfrac{v_k}{2}.$$

The following lemma ensures that the pay-off matrices $A$ and $B$ satisfy the assumptions in Theorems 2.1 and 2.2. The proof is in the Appendix.

LEMMA 5.1. *For the war of attrition with pay-off matrix given by* (5.1) *or* (5.3) *there exists a unique ESS, and the pay-off matrix is conditionally negative definite.*

The next theorem is a basic persistence result for the replicator dynamics, saying that with probability close to 1 the maximum effort strategy, that is, strategy $n$, will not die out. For $j = 0, \ldots, n$, let $\mathbf{e}_j$ denote the $(j+1)$st column of the unit matrix of order $n+1$.

THEOREM 5.1. *Let $X(t)$ be given by the stochastic replicator dynamics* (5.2) *with initial state $\mathbf{x} \in \Delta$. Let $\varepsilon > 0$. Then there exists $\sigma^* = \sigma^*(\varepsilon) > 0$ such that*

$$P_{\mathbf{x}}\left\{\limsup_{t \to \infty} X_n(t) > 0\right\} \geq 1 - \varepsilon,$$

*provided that $\sigma_0, \ldots, \sigma_n < \sigma^*$.*

PROOF. Let $\mathbf{p}$ be the ESS for $A$. It will first be shown that $p_n > 0$. Suppose that this is not the case, so that $m := \max\{j : p_j > 0\} < n$. Then $\mathbf{e}_{m+1}^T A \mathbf{p} \leq \mathbf{e}_m^T A \mathbf{p}$, see Hofbauer and Sigmund [(1998), page 63]. As $p_j = 0$ for $j > m$,

$$\mathbf{e}_m^T A \mathbf{p} = \sum_{j=0}^{m-1}(v_j - c_j)p_j + \left(\frac{v_m}{2} - c_m\right)p_m,$$

$$\mathbf{e}_{m+1}^T A \mathbf{p} = \sum_{j=0}^{m}(v_j - c_j)p_j.$$

Thus

$$0 \leq \mathbf{e}_m^T A \mathbf{p} - \mathbf{e}_{m+1}^T A \mathbf{p} = -\frac{v_m}{2}p_m < 0;$$

a contradiction. Hence $p_n > 0$.

By Theorem 2.1 and Lemma 5.1, for all $t > 0$,

$$E_{\mathbf{x}}\frac{1}{t}\int_0^t |X_n(s) - p_n|^2\,ds \leq \frac{1}{|\lambda_2|}\left\{\frac{d(\mathbf{x}, \mathbf{p})}{t} + \frac{1}{2}\sum_{j=1}^n p_j \sigma_j^2\right\},$$

where $\lambda_2 \neq 0$ depends only on $A$. Choose $t_0 > 0$ and $\sigma^* > 0$ such that

$$\frac{d(\mathbf{x}, \mathbf{p})}{|\lambda_2|t_0} < \frac{p_n^2\varepsilon}{16}, \qquad \frac{(\sigma^*)^2}{|\lambda_2|} < \frac{p_n^2\varepsilon}{8}.$$

Thus if $\sigma_0, \ldots, \sigma_n < \sigma^*$, then

$$E_{\mathbf{x}}\frac{1}{t}\int_0^t |X_n(s) - p_n|^2\,ds \leq \frac{p_n^2\varepsilon}{8} \qquad \text{for all } t \geq t_0.$$

Now consider the increasing sequence of events

$$F_\mu = \left\{ X_n(s) \leq \frac{p_n}{2} \text{ for all } s \geq t_0 + \mu \right\}, \qquad \mu = 1, 2, \ldots.$$

For every $\mu$,

$$\frac{p_n^2 \varepsilon}{8} \geq \frac{1}{2(t_0 + \mu)} \int_{F_\mu} \int_{t_0+\mu}^{2(t_0+\mu)} |X_n(s) - p_n|^2 \, ds \, dP_\mathbf{x} \geq \frac{p_n^2}{8} P_\mathbf{x}(F_\mu),$$

so that $P_\mathbf{x}(\bigcup_{\mu=1}^\infty F_\mu) \leq \varepsilon$. Hence

$$P_\mathbf{x}\left\{ \limsup_{t \to \infty} X_n(t) > 0 \right\} \geq P_\mathbf{x}\left( \bigcap_{\mu=1}^\infty F_\mu^C \right) \geq 1 - \varepsilon. \qquad \square$$

In general the ESS for a discrete war of attrition may have a fairly complicated structure. Cressman [(2003), Section 7.4] describes a broad approach to calculating ESSs for these games using backward induction. Whittaker (1996) has recently solved several closely related resource allocation problems based on a multiple trial war of attrition. The following theorem gives an explicit expression for the ESS in the case where the $\rho_k$ and the $v_k$ are constant and $c_j = j$ for all strategies $j$. Combining Lemma 5.1 and Theorems 2.1, 2.2 and 5.2, one obtains a fairly complete picture of the long-run behavior of the stochastic replicator dynamics when the conflicts are modeled by a war of attrition.

Let $U_m(x)$ denote the $m$th Chebyshev polynomial of the second kind, and let $U_{-1}(x) \equiv 0$. Let $i = \sqrt{-1}$.

THEOREM 5.2.  *Consider the war of attrition with pay-off matrix $B = (b_{jk})$, $j, k = 0, \ldots, n$, where*

$$(5.4) \qquad b_{jk} = \begin{cases} v - k, & j > k, \\ \dfrac{v}{2} - k - \rho, & j = k, \\ -j, & j < k, \end{cases}$$

*and $0 \leq \rho < \frac{1}{2}v$. The unique ESS $\mathbf{p}$ is given as follows. If the reward $v$ is so large that $v \geq 2n + 2\rho$, then*

$$p_0 = \cdots = p_{n-1} = 0, \qquad p_n = 1.$$

*Otherwise there is a unique index $s \in \{0, \ldots, n-1\}$ such that*

$$(5.5) \qquad n - 1 + \rho \leq \frac{v}{2} + s < n + \rho,$$

*and*

$$p_k = \frac{1}{c}\left(-\frac{v}{2} - \rho\right)^k$$

(5.6)
$$\times \left\{ u_{s-k+1} + \left(s + 1 - n + \frac{v}{2} + \rho\right) u_{s-k} \right.$$
$$\left. + (s+1-n)\left(\frac{v}{2} + \rho\right) u_{s-k-1} \right\}, \quad 0 \le k \le s,$$

(5.7) $\quad p_k = 0, \quad s+1 \le k \le n-1,$

$$p_n = \frac{1}{c}\left(-\frac{v}{2} - \rho\right)^{s+1},$$

*where*

$$u_k = (-i\gamma)^k U_k\left(-\frac{i(2\rho + 1)}{2\gamma}\right), \quad \gamma = \sqrt{\frac{v^2}{4} - \rho^2}$$

*and*

$$c = -u_{s+2} + (n - s - 1 - 2\rho)u_{s+1} + \{2\rho(n - s - 1) + \gamma^2\}u_s$$
$$- (n - s - 1)\gamma^2 u_{s-1}.$$

The proof of this theorem requires some auxiliary results, which are proved in the Appendix.

It was shown in the proof of Theorem 5.1 that for the war of attrition (5.3), strategy $n$ is always contained in the support of the ESS. The next lemma states that for $j < n$, strategy $j$ can belong to the support only if the corresponding cost $c_j$ is below a certain threshold. This is the discrete analogue of Theorem 7 of Bishop and Cannings (1978). The lemma explains in particular the choice of $s$ in Theorem 5.2.

LEMMA 5.2. *Let* **p** *be the ESS for the war of attrition with pay-off matrix* $B = (b_{jk})$, $j, k = 0, \ldots, n$, *described by* (5.3). *If* $j < n$ *and* $c_j \ge c_n + \rho_n - \frac{1}{2}v_n$, *then* $p_j = 0$.

The next two lemmas give explicit formulas for determinants related to the pay-off matrix of a war of attrition. Let $J_k$ denote the $k \times k$ matrix with all entries equal to 1 and let $\mathbf{1}_k$ denote the $k \times 1$ vector all of whose entries are 1.

LEMMA 5.3. *Let* $B \in \mathbb{R}^{(n+1)\times(n+1)}$ *be given by* (5.3). *For* $k = 0, \ldots, n$, *let* $B^{(k)}$ *denote the matrix obtained from* $B$ *by replacing column* $k$ *with the* *vector* $\mathbf{1}_{n+1}$. *Then*

$$\det B^{(n)} = \prod_{j=0}^{n-1}\left(-\frac{v_j}{2} - \rho_j\right)$$

*and, for* $k = 0, \ldots, n-1$,

$$\det B^{(k)} = \det(\widetilde{B}_{n-k} + c_k J_{n-k}) \prod_{j=0}^{k-1}\left(-\frac{v_j}{2} - \rho_j\right),$$

*where* $\widetilde{B}_{n-k}$ *is the* $(n-k) \times (n-k)$ *principal submatrix of* $B$ *situated in the bottom right-hand corner.*

LEMMA 5.4. *Let* $B \in \mathbb{R}^{(n+1)\times(n+1)}$ *be given by* (5.4) *with* $0 \leq \rho < \frac{1}{2}v$. *Set* $\gamma = \sqrt{\frac{1}{4}v^2 - \rho^2}$. *Then*

$$\det B = \left(\frac{v}{2} - \rho\right)(-i\gamma)^{n-1}$$
$$\times \left\{-i\gamma U_n\left(-\frac{i(2\rho+1)}{2\gamma}\right) + \left(\frac{v}{2} + \rho\right)U_{n-1}\left(-\frac{i(2\rho+1)}{2\gamma}\right)\right\}.$$

PROOF OF THEOREM 5.2. Let $B$ be given by (5.4). Suppose first that $v \geq 2n + 2\rho$. Then $j \geq n + \rho - \frac{1}{2}v$ for every $j = 0, \ldots, n-1$. Thus if $\mathbf{p}$ is the ESS, then, by Lemma 5.2, $p_j = 0$ for $j = 0, \ldots, n-1$, so that $\mathbf{p} = \mathbf{e}_n$.

Suppose next that $2\rho < v < 2n + 2\rho$. Define $s \in \{0, \ldots, n-1\}$ by (5.5) and define $\mathbf{p} \in \mathbb{R}^{n+1}$ by (5.6), (5.7) and (5.8). It will be shown that

(5.8) $\quad\quad \{B\mathbf{p}\}_j = \{B\mathbf{p}\}_n \quad$ if $0 \leq j \leq s$,

(5.9) $\quad\quad \{B\mathbf{p}\}_j \leq \{B\mathbf{p}\}_n \quad$ if $s+1 \leq j \leq n-1$.

It will also be shown that $\mathbf{p} \in \overline{\Delta}$. Since, by Lemma 5.1, any principal submatrix of $B$ is conditionally negative definite, it will then follow from Haigh's (1975) theorem that $\mathbf{p}$ is the ESS.

For $m = 0, \ldots, n$ set

$$B_m = \begin{bmatrix} b_{00} & \cdots & b_{0m} \\ \vdots & & \vdots \\ b_{m0} & \cdots & b_{mm} \end{bmatrix}.$$

According to Lemma 5.4,

(5.10) $\quad\quad \det B_m = \left(\frac{v}{2} - \rho\right)\left\{u_m + \left(\frac{v}{2} + \rho\right)u_{m-1}\right\}.$

Set
$$\overline{B} = \begin{bmatrix} b_{00} & \cdots & b_{0s} & b_{0n} \\ \vdots & & \vdots & \vdots \\ b_{s0} & \cdots & b_{ss} & b_{sn} \\ b_{n0} & \cdots & b_{ns} & b_{nn} \end{bmatrix} \quad \text{and} \quad \overline{\mathbf{p}} = \begin{pmatrix} p_0 \\ \vdots \\ p_s \\ p_n \end{pmatrix}.$$

For $k = 0, \ldots, s+1$ let $\overline{B}^{(k)}$ denote the matrix obtained from $\overline{B}$ by replacing column $k$ with the vector $\mathbf{1}_{s+2}$. As $p_k = 0$ for $s+1 \leq k \leq n-1$,

(5.11) $\quad \{\overline{B}\overline{\mathbf{p}}\}_j = \{B\mathbf{p}\}_j, \quad 0 \leq j \leq s, \quad \{\overline{B}\overline{\mathbf{p}}\}_{s+1} = \{B\mathbf{p}\}_n.$

The matrix $\overline{B}$ is again of the form (5.3). It follows from Lemma 5.3 that for every $k \leq s-1$,

$$\det \overline{B}^{(k)} = \left(-\frac{v}{2} - \rho\right)^k \det \begin{bmatrix} b_{k+1,k+1} + k & \cdots & b_{k+1,s} + k & b_{k+1,n} + k \\ \vdots & & \vdots & \vdots \\ b_{s,k+1} + k & \cdots & b_{s,s} + k & b_{s,n} + k \\ b_{n,k+1} + k & \cdots & b_{n,s} + k & b_{n,n} + k \end{bmatrix}$$

$$= \left(-\frac{v}{2} - \rho\right)^k$$

$$\times \det \begin{bmatrix} \frac{v}{2} - \rho - 1 & -1 & \cdots & -1 & -1 \\ v-1 & \frac{v}{2} - \rho - 2 & \cdots & -2 & -2 \\ \vdots & \vdots & & \vdots & \vdots \\ v-1 & v-2 & \cdots & \frac{v}{2} - \rho - (s-k) & -(s-k) \\ v-1 & v-2 & \cdots & v - (s-k) & \frac{v}{2} - \rho - (n-k) \end{bmatrix}.$$

Denote the matrix in the previous line by $Q$. To calculate $\det Q$ augment $Q$ from the left by the $(s-k+1) \times 1$ vector $(v, \ldots, v)^T$ and put on top of the matrix thus obtained the $1 \times (s-k+2)$ vector $(\frac{1}{2}v - \rho, 0, \ldots, 0)$. This gives a matrix which is equal to $B_{s-k+1}$ except for the element in the lower right-hand corner, which is $\frac{1}{2}v - \rho - n + k$, while the corresponding element of $B_{s-k+1}$ is $\frac{1}{2}v - \rho - s + k - 1$. Hence

$$\det Q = \left(\frac{v}{2} - \rho\right)^{-1} \det \begin{bmatrix} \frac{v}{2} - \rho & \mathbf{0}_{1 \times (s-k+1)} \\ v\mathbf{1}_{s-k+1} & Q \end{bmatrix}$$

$$= \left(\frac{v}{2} - \rho\right)^{-1} \{\det B_{s-k+1} - (n-s-1)\det B_{s-k}\},$$

and so, by (5.10),

$$\det \overline{B}^{(k)} = \left(-\frac{v}{2} - \rho\right)^k \left[u_{s-k+1} + \left(\frac{v}{2} + \rho\right) u_{s-k}\right.$$

$$-(n-s-1)\left\{u_{s-k} + \left(\frac{v}{2}+\rho\right)u_{s-k-1}\right\}\right]$$

$$= cp_k,$$

provided $k \leq s-1$. Similarly, $\det \overline{B}^{(s)} = cp_s$, and it follows directly from Lemma 5.3 that

$$\det \overline{B}^{(s+1)} = \left(-\frac{v}{2}-\rho\right)^{s+1} = cp_n.$$

Thus

(5.12) $$c\bar{\mathbf{p}} = \begin{pmatrix} \det \overline{B}^{(0)} \\ \vdots \\ \det \overline{B}^{(s+1)} \end{pmatrix}.$$

Let $\operatorname{adj} \overline{B}$ denote the adjugate matrix of $\overline{B}$. It is readily verified that

(5.13) $$(\operatorname{adj} \overline{B})\mathbf{1}_{s+2} = \begin{pmatrix} \det \overline{B}^{(0)} \\ \vdots \\ \det \overline{B}^{(s+1)} \end{pmatrix},$$

and therefore

$$\overline{B}\bar{\mathbf{p}} = \frac{1}{c}\overline{B}(\operatorname{adj} \overline{B})\mathbf{1}_{s+2} = \frac{\det \overline{B}}{c}\mathbf{1}_{s+2}.$$

In view of (5.11) this proves the first claim (5.8).

If $s < j < n$, then, by (5.5),

$$b_{jn} = -j \leq -s-1 \leq \frac{v}{2} - n - \rho = b_{nn},$$

and so

$$\{B\mathbf{p}\}_j = \sum_{k=0}^{s} b_{jk}p_k + b_{jn}p_n = \sum_{k=0}^{s} b_{nk}p_k + b_{jn}p_n \leq \sum_{k=0}^{s} b_{nk}p_k + b_{nn}p_n = \{B\mathbf{p}\}_n,$$

proving the second claim (5.9).

Finally, to show that $\mathbf{p} \in \overline{\Delta}$ set

$$t_j = u_j + \left(s+1-n+\frac{v}{2}+\rho\right)u_{j-1} + (s+1-n)\left(\frac{v}{2}+\rho\right)u_{j-2}, \qquad j=1,2,\ldots,$$

so that the first $s+1$ entries of $\mathbf{p}$ can be written as

$$p_k = \frac{1}{c}\left(-\frac{v}{2}-\rho\right)^k t_{s-k+1}, \qquad 0 \leq k \leq s.$$

In view of (5.5),

$$t_1 = \frac{v}{2} - \rho - n + s < 0,$$

$$t_2 = \left(\frac{v}{2} + \rho\right)\left(\frac{v}{2} - \rho - n + s + 1\right) - (2\rho + 1)t_1 > 0.$$

Using the recurrence relation for the Chebyshev polynomials [Szegö (1975), equation (4.7.17), page 81] one may verify that

$$u_k = -(2\rho + 1)u_{k-1} + \gamma^2 u_{k-2}, \qquad k \geq 1.$$

Therefore,

$$t_k = -(2\rho + 1)t_{k-1} + \gamma^2 t_{k-2}, \qquad k \geq 3.$$

It now follows by induction that $(-1)^k t_k > 0$ for all $k$. A short calculation shows that $c = -t_{s+2} + (\frac{1}{2}v - \rho)t_{s+1}$, so that $(-1)^{s+1}c > 0$. (In particular, $c \neq 0$, which has hitherto been tacitly assumed.) It is thus obvious that $p_k \geq 0$ for $0 \leq k \leq n$.

To verify that $\sum_{k=0}^n p_k = 1$ note that by (5.7), (5.12) and (5.13),

$$\sum_{k=0}^n p_k = p_n + \sum_{k=0}^s p_k = \frac{1}{c}\mathbf{1}_{s+2}^T(\operatorname{adj}\overline{B})\mathbf{1}_{s+2}.$$

By a well-known determinantal formula for partitioned matrices [Gantmacher (1959), page 46],

$$\mathbf{1}_{s+2}^T(\operatorname{adj}\overline{B})\mathbf{1}_{s+2} = \det\overline{B} - \det(\overline{B} - \mathbf{1}_{s+2}\mathbf{1}_{s+2}^T).$$

Observing that, by (5.10),

$$\det\overline{B} = \det B_{s+1} - (n - s - 1)\det B_s$$
$$= \left(\frac{v}{2} - \rho\right)\left[u_{s+1} + \left(\frac{v}{2} + \rho\right)u_s - (n - s - 1)\left\{u_s + \left(\frac{v}{2} + \rho\right)u_{s-1}\right\}\right]$$

and

$$\det(\overline{B} - \mathbf{1}_{s+2}\mathbf{1}_{s+2}^T) = \left(\frac{v}{2} - \rho\right)^{-1}\{\det B_{s+2} - (n - s - 1)\det B_{s+1}\}$$
$$= u_{s+2} + \left(\frac{v}{2} + \rho\right)u_{s+1} - (n - s - 1)\left\{u_{s+1} + \left(\frac{v}{2} + \rho\right)u_s\right\},$$

one obtains that $\mathbf{1}_{s+2}^T(\operatorname{adj}\overline{B})\mathbf{1}_{s+2} = c$, so that $\sum_{k=0}^n p_k = 1$.  $\square$

## APPENDIX

PROOF OF LEMMA 5.1. As $\mathbf{q}^T B \mathbf{q} = \mathbf{q}^T A \mathbf{q} - \sum_{j=0}^{n} q_j^2 \rho_j \leq \mathbf{q}^T A \mathbf{q}$ for all $\mathbf{q} \in \mathbb{R}^{n+1}$, it suffices to show that $A$ is conditionally negative definite. Define the $n \times n$ matrix $D$ by $d_{jk} = \frac{1}{2}(b_{jk} + b_{kj}) - b_{j0} - b_{0k} + b_{00}$, $j,k = 1,\ldots,n$. Thus

$$d_{jk} = v_{\min\{j,k\}} - 2c_{\min\{j,k\}} - v_0 + 2c_0.$$

For $k = 1,\ldots,n$ let $\mathbf{f}_k$ be the $n \times 1$ vector whose first $k-1$ entries are 0 and whose remaining entries are 1. Then $D$ can be written as

$$D = \sum_{k=1}^{n} (v_k - v_{k-1} - 2(c_k - c_{k-1})) \mathbf{f}_k \mathbf{f}_k^T,$$

showing that $D$ is negative definite. This implies that $A$ is conditionally negative definite; see Haigh (1975).

The existence of an ESS now follows, since in a game with conditionally negative definite pay-off matrix, a Nash equilibrium, which always exists, must be an ESS. To prove uniqueness suppose $\mathbf{p}$ and $\mathbf{q}$ are ESSs. Then $\mathbf{p}^T B \mathbf{p} \geq \mathbf{q}^T B \mathbf{p}$ and $\mathbf{q}^T B \mathbf{q} \geq \mathbf{p}^T B \mathbf{q}$, so that $(\mathbf{p} - \mathbf{q})^T B (\mathbf{p} - \mathbf{q}) \geq 0$. Since $B$ is conditionally negative definite, this implies $\mathbf{p} = \mathbf{q}$. □

PROOF OF LEMMA 5.2. The assertion is obviously true if $p_j = 0$ for all $j < n$, so assume $p_j > 0$ for some $j < n$. Let $m := \max\{j : j < n,\ p_j > 0\}$. Then $p_m > 0$ and, as in the proof of Theorem 5.1, $p_n > 0$, so that $\mathbf{e}_m^T B \mathbf{p} = \mathbf{e}_n^T B \mathbf{p}$, see Hofbauer and Sigmund [(1998), page 63]. Since $p_j = 0$ if $m + 1 \leq j \leq n - 1$,

$$\mathbf{e}_m^T B \mathbf{p} = \sum_{k=0}^{m-1} (v_k - c_k) p_k + \left(\frac{v_m}{2} - c_m - \rho_m\right) p_m - c_m p_n,$$

$$\mathbf{e}_n^T B \mathbf{p} = \sum_{k=0}^{m} (v_k - c_k) p_k + \left(\frac{v_n}{2} - c_n - \rho_n\right) p_n.$$

Thus

$$0 = \mathbf{e}_n^T B \mathbf{p} - \mathbf{e}_m^T B \mathbf{p} = \left(\frac{v_m}{2} + \rho_m\right) p_m + \left(\frac{v_n}{2} - c_n - \rho_n + c_m\right) p_n,$$

and it follows that $v_n/2 - c_n - \rho_n + c_m < 0$. That is, $c_n + \rho_n - v_n/2 > c_m$. Now suppose that $j < n$ and $c_j \geq c_n + \rho_n - v_n/2$. Then $c_j > c_m$, and since the sequence $(c_\mu)$ is increasing, $j > m$. Thus $p_j = 0$. □

PROOF OF LEMMA 5.3. Suppose $1 \leq k \leq n-1$. For $j = 0,\ldots,k-1$, add $v_j - c_j$ times column $k$ of $B^{(k)}$, that is, the vector $\mathbf{1}_{n+1}$, to column $j$.

For $j = k+1, \ldots, n$, add $c_k$ times column $k$ to column $j$. The matrix thus obtained can be partitioned as

$$\begin{bmatrix} D & \mathbf{1}_k & * \\ \mathbf{0} & 1 & \mathbf{0} \\ 0 & \mathbf{1}_{n-k} & \widetilde{B}_{n-k} + c_k J_{n-k} \end{bmatrix},$$

where $D$ is a $k \times k$ upper triangular matrix with diagonal elements $-v_0/2 - \rho_0, \ldots, -v_{k-1}/2 - \rho_{k-1}$. The assertion is now obvious. The proof is similar for $k=0$ and $k=n$. $\square$

LEMMA A.1. *Let $\gamma_1, \gamma_2, x \in \mathbb{R}$, $\gamma_1, \gamma_2 > 0$. The determinant of the $n \times n$ tridiagonal matrix*

$$D_n(x) = \begin{bmatrix} x & \gamma_1 & 0 & 0 & \ldots & 0 & 0 & 0 \\ -\gamma_2 & x & \gamma_1 & 0 & \ldots & 0 & 0 & 0 \\ 0 & -\gamma_2 & x & \gamma_1 & \ldots & 0 & 0 & 0 \\ \ldots & \ldots & \ldots & \ldots & \ldots & \ldots & \ldots & \ldots \\ 0 & 0 & 0 & 0 & \ldots & -\gamma_2 & x & \gamma_1 \\ 0 & 0 & 0 & 0 & \ldots & 0 & -\gamma_2 & x \end{bmatrix}$$

*is given by*

(A.1) $$\det D_n(x) = i^n (\gamma_1 \gamma_2)^{n/2} U_n\left(-\frac{ix}{2\sqrt{\gamma_1 \gamma_2}}\right).$$

PROOF. Note first that

$$\det D_1(x) = x, \qquad \det D_2(x) = x^2 + \gamma_1 \gamma_2.$$

Expanding $\det D_n(x)$ along the last column, one obtains that, for $n > 2$,

$$\det D_n(x) = x \det D_{n-1}(x) + \gamma_1 \gamma_2 \det D_{n-2}(x).$$

Denote the expression on the right-hand side of (A.1) by $h_n(x)$. Then

$$h_1(x) = x, \qquad h_2(x) = x^2 + \gamma_1 \gamma_2,$$

and, by the recurrence formula for the Chebyshev polynomials [Szegö (1975), (4.7.17), page 81], for $n > 2$,

$$h_n(x) = i^{n-1}(\gamma_1\gamma_2)^{(n-1)/2} x U_{n-1}\left(-\frac{ix}{2\sqrt{\gamma_1\gamma_2}}\right) - i^n(\gamma_1\gamma_2)^{n/2} U_{n-2}\left(-\frac{ix}{2\sqrt{\gamma_1\gamma_2}}\right)$$
$$= x h_{n-1}(x) + \gamma_1 \gamma_2 h_{n-2}(x).$$

Now the assertion follows by induction. $\square$

PROOF OF LEMMA 5.4. Define the $n \times n$ matrices $F$ and $G$ by

$$F = \begin{bmatrix} \frac{v}{2} - \rho & & & & \\ v & \frac{v}{2} - \rho & & & \\ v & v & \frac{v}{2} - \rho & & \\ \vdots & \vdots & & & \\ v & v & v & \cdots & v & \frac{v}{2} - \rho \end{bmatrix},$$

$$G = \begin{bmatrix} -1 & -1 & -1 & \cdots & -1 \\ -1 & -2 & -2 & \cdots & -2 \\ -1 & -2 & -3 & \cdots & -3 \\ \vdots & \vdots & \vdots & & \vdots \\ -1 & -2 & -3 & \cdots & -n \end{bmatrix}.$$

Let $I$ denote the $n \times n$ unit matrix. Then

$$\det B = \left(\frac{v}{2} - \rho\right) \det(F + G) = \left(\frac{v}{2} - \rho\right) \det G \det(G^{-1}F + I).$$

It is easily seen that $\det G = (-1)^n$. Moreover,

$$G^{-1} = \begin{bmatrix} -2 & 1 & & & & \\ 1 & -2 & 1 & & & \\ & 1 & -2 & 1 & & \\ & & \ddots & \ddots & \ddots & \\ & & & 1 & -2 & 1 \\ & & & & 1 & -1 \end{bmatrix},$$

so that

$$G^{-1}F + I = \begin{bmatrix} 2\rho+1 & \frac{v}{2} - \rho & & & & \\ -\frac{v}{2} - \rho & 2\rho+1 & \frac{v}{2} - \rho & & & \\ & -\frac{v}{2} - \rho & 2\rho+1 & \frac{v}{2} - \rho & & \\ & & \ddots & & \ddots & \\ & & & -\frac{v}{2} - \rho & 2\rho+1 & \frac{v}{2} - \rho \\ & & & & -\frac{v}{2} - \rho & -\frac{v}{2} + \rho + 1 \end{bmatrix}.$$

It now follows by Lemma A.1 that

$$\det B = \left(\frac{v}{2} - \rho\right)(-1)^n \left\{ i^n \gamma^n U_n\left(-\frac{i(2\rho+1)}{2\gamma}\right) \right.$$
$$\left. - \left(\frac{v}{2} + \rho\right) i^{n-1} \gamma^{n-1} U_{n-1}\left(-\frac{i(2\rho+1)}{2\gamma}\right) \right\}$$

$$= (-i\gamma)^n \left\{ \left(\frac{v}{2} - \rho\right) U_n\left(-\frac{i(2\rho+1)}{2\gamma}\right) + i\gamma U_{n-1}\left(-\frac{i(2\rho+1)}{2\gamma}\right) \right\}. \quad \square$$

Institut für Statistik
Aachen University
D-52056 Aachen
Germany
e-mail: imhof@stochastik.rwth-aachen.de